\documentclass[12pt,dvips]{article}
\usepackage{amsmath,amsfonts,amssymb,amsthm,graphicx,verbatim,subfig}
\usepackage[all]{xy}

\ifx\pdftexversion\undefined

\usepackage[a4paper,colorlinks,link
=black,filecolor=black,citecolor=black,urlcolor=black,pdfstartview=FitH]{hyperref}
\else

\usepackage[a4paper,colorlinks,linkcolor=black,filecolor=black,citecolor=black,urlcolor=black,pdfstartview=FitH]{hyperref}
\fi

\font\sixbb=msbm6
\font\eightbb=msbm8
\font\twelvebb=msbm10 scaled 1095
\newfam\bbfam
\textfont\bbfam=\twelvebb \scriptfont\bbfam=\eightbb
                           \scriptscriptfont\bbfam=\sixbb

\newcommand{\Rea}{{\mathbb R}}

\newcommand{\Int}{{\mathbb Z}}

\newcommand{\Nat}{{\mathbb N}}
\newcommand{\FF}{{\mathbb F}}

\newtheorem{theorem}{\bf Theorem}[section]
\newtheorem{claim}[theorem]{\bf Claim}
\newtheorem{conjecture}[theorem]{\bf Conjecture}
\newtheorem{proposition}[theorem]{\bf Proposition}
\newtheorem{corollary}[theorem]{\bf Corollary}
\newtheorem{definition}[theorem]{\bf Definition}
\newcommand{\enp}{\begin{flushright} $\Box$ \end{flushright}}
\newcommand{\beq}[0]{\begin{equation}}
\newcommand{\enq}[0]{\end{equation}}

\newcommand{\cb}{{\cal B}}

\newcommand{\thh}{\tilde{H}}
\newcommand{\supp}{{\rm supp}}
\newcommand{\cf}{{\cal F}}

\newcommand{\stab}{{\rm stab}}
\newcommand{\bst}{B_{s,\tau}}
\newcommand{\isig}{\iota_{s}}

\newcommand{\zhat}{\widehat{0}}
\newcommand{\ohat}{\widehat{1}}
\newcommand{\ccs}{{\mathbb S}}

\newcommand{\scst}{\supp(c_{s,\tau})}
\newcommand{\tlam}{\tilde{\lambda}}
\newcommand{\ccc}{{\cal C}}
\newcommand{\uv}{\underline{v}}
\newcommand{\ctil}{\tilde{c}}

\title{Expansion of Building-Like Complexes}
\begin{document}
\author{Alexander Lubotzky\thanks{Institute of Mathematics, Hebrew University, Jerusalem 91904, Israel. e-mail: alexlub@math.huji.ac.il~. Supported by ERC, ISF and NSF grants.} \and Roy Meshulam\thanks{Department of Mathematics,
Technion, Haifa 32000, Israel. e-mail:
meshulam@math.technion.ac.il~. Supported by an ISF grant.} \and Shahar Mozes\thanks{Institute of Mathematics, Hebrew University, Jerusalem 91904, Israel. e-mail: mozes@math.huji.ac.il~. Supported by ISF and BSF grants}}
\maketitle
\pagestyle{plain}
\begin{abstract}
Following Gromov, the coboundary expansion of building-like complexes is studied.
In particular, it is shown that for any $n \geq 1$, there exists a constant $\epsilon(n)>0$ such that for any $0 \leq k <n$ the $k$-th coboundary expansion constant of any $n$-dimensional spherical building is at least $\epsilon(n)$.

\end{abstract}

\section{Introduction}
\label{intro}
Expander graphs have been a focus of intensive research in the last four decades, with many applications in combinatorics and computer science as well as pure mathematics
(see \cite{HLW06,Lu1,Lu2}). In recent years a high dimensional theory is emerging. There are several ways to extend the definition of expanders from graphs to simplicial complexes (see \cite{Lu3} for a survey). Here we will be concerned with the notion of "coboundary expansion" that came up independently in the work of Linial, Meshulam and Wallach \cite{LM06,MW09} on homological connectivity of random complexes and in Gromov's work \cite{G10} on the topological overlap property. For an application of coboundary expansion to property testing see
\cite{KL13}. The rich theory of expander graphs hints that high dimensional expanders can also be useful. The goal of this paper is to show, following Gromov \cite{G10}, that spherical buildings, and more generally, "building-like complexes" (defined precisely below), are expanders.

We proceed with the formal definitions.
Let $X$ be a finite $n$-dimensional pure simplicial complex. For $k \geq 0$, let $X^{(k)}$ denote the $k$-dimensional skeleton of $X$ and let $X(k)$ be the family of $k$-dimensional faces of $X$, $f_k(X)=|X(k)|$.
Define a positive weight function $w=w_X$ on the simplices of $X$ as follows. For
$\sigma \in X(k)$, let $c(\sigma)=|\{\eta \in X(n): \sigma \subset \eta\}|$
and let $$w(\sigma)=\frac{c(\sigma)}{\binom{n+1}{k+1} f_{n}(X)}.$$
Note that $\sum_{\sigma \in X(k)} w(\sigma)=1$ and
if $\sigma \in X(k)$ then $$\sum_{\{\tau \in X(k+1): \sigma \subset \tau\}} w(\tau)=(k+2)w(\sigma).$$
All homology and cohomology groups referred to in the sequel are with $\FF_2$ coefficients.
Let $C_k(X)$ be the space of $\FF_2$-valued $k$-chains of $X$ with the boundary map $\partial_k:C_k(X) \rightarrow C_{k-1}(X)$. Let $C^k(X)$ denote the space of  $\FF_2$-valued $k$-cochains of $X$ with the coboundary map $d_k:C^k \rightarrow C^{k+1}$.
As usual, the spaces of $k$-cycles and $k$-cocycles are denoted by $Z_k(X)$ and $Z^k(X)$ and the spaces of $k$-boundaries and $k$-coboundaries are denoted by $B_k(X)$ and $B^k(X)$.
Reduced $k$-dimensional homology and cohomology will be denoted by $\thh_k(X)$ and $\thh^k(X)$.
For $\phi \in C^k(X)$, let $[\phi]$ denote the image of $\phi$ in
$C^k(X)/B^k(X)$.
Let $$\|\phi\|=\sum_{\{\sigma \in X(k): \phi(\sigma) \neq 0\}} w(\sigma)$$
and  $$\|[\phi]\|=\min\{\|\phi+d_{k-1}\psi\|: \psi \in C^{k-1}(X)\}.$$
\begin{definition}
\label{hexp}
The {\it $k$-th coboundary expansion constant} of $X$ is
$$
h_k(X)=\min \left\{\frac{\|d_k\phi \|}{\|[\phi]\|}: \phi \in C^k(X)-B^k(X)\right\}.
$$
\end{definition}

\noindent
{\bf Remarks:}\\
1. Note that $h_k(X)=0$ iff $\tilde{H}^k(X;\FF_2) \neq 0$.
\\
2. Let $\Delta_n$ denote the $n$-simplex and let $0 \leq k \leq n-1$.
In \cite{MW09,G10} it was shown that the $k$-th coboundary expansion of $\Delta_n$ satisfies
\begin{equation}
\label{expsim}
h_{k}(\Delta_n) \geq \frac{n+1}{n-k}
\end{equation}
with equality when $n+1$ is divisible by $k+2$.
\\
3. Let $k<n$ and let $\sigma \in X(k)$ be a $k$-simplex of minimal weight. Then $\|[1_{\sigma}]\|=w(\sigma)$ and therefore
\begin{equation}
\label{simp1}
h_k(X) \leq \frac{\|d_k 1_{\sigma} \|}{\|[1_{\sigma}]\|} =
\frac{\sum \{w(\tau):\tau \in X(k+1), \sigma \subset \tau\}}{w(\sigma)}=k+2.
\end{equation}
Equality in (\ref{simp1}) is attained for $X=\Delta_n$ and $k=n-1$.
\\
4. The normalization we use for the norm in $C^k(X)$ and hence for the definition of $h_k(X)$ takes into account the possibility that
the $k$-faces of $X$ may not all have the same degrees. This is particularly relevant for spherical buildings - see the example
following Corollary \ref{bilexp}.

In this note we are concerned with the expansion of certain building-like complexes.
Let $G$ be a subgroup of $Aut(X)$ and let $S$ be a finite $G$-set.
For $0 \leq k \leq n-1$, let $\cf_k=S \times X(k)$ with a $G$-action given by
$g(s,\tau)=(gs, g\tau)$.
Let $$\cb=\{\bst: -1 \leq k <n, (s,\tau) \in \cf_k\}$$ be a family of subcomplexes of $X$ such that
$\tau \in B_{s,\tau} \subset B_{s,\tau'}$ for all $s \in S$ and
$\tau \subset \tau' \in X^{(n-1)}$.
\begin{definition}
\label{blcom}
A {\it building-like complex} is a $4$-tuple $(X,S,G,\cb)$ as above with the following properties:
\begin{enumerate}
\item[(C1)]
$G$ is transitive on $X(n)$.
\item[(C2)]
$g\bst=B_{gs,g\tau}$ for all $g \in G$ and $(s,\tau) \in S \times X^{(n-1)}$.
\item[(C3)]
$\thh_i(\bst)=0$ for all $(s,\tau) \in \cf_k$ and $-1 \leq i \leq k<n$.
\end{enumerate}
\end{definition}
Examples of building-like complexes include basis-transitive matroid complexes and spherical buildings - see Section \ref{blike}.
Following Gromov \cite{G10}, we give a lower bound on the expansion of building-like complexes.
For a simplex $\eta \in X$, let $G\eta$ denote the orbit of $\eta$ under $G$.
For $0 \leq k \leq n-1$, let
$$a_k=a_k(X,S,G,\cb)=\max\{|G\eta \cap \bst(k+1)|: \eta \in X(k+1), (s,\tau) \in \cf_{k}\}.$$
\begin{theorem}
\label{gromov}
Let $(X,S,G,\cb)$ be an $n$-dimensional building-like complex.  Then
for $0 \leq k \leq n-1$,
\begin{equation}
\label{bound}
h_{k}(X) \geq \left(\binom{n+1}{k+2}a_{k}\right)^{-1}.
\end{equation}
\end{theorem}

The proof of Theorem \ref{gromov} is given in Section \ref{prf}.
In Section \ref{blike} we use Theorem \ref{gromov} to derive expansion bounds for basis-transitive matroid complexes and for spherical buildings.
In Section \ref{appl} we discuss applications to topological overlapping and to property testing.
We conclude in Section \ref{conc} with some questions and comments.

\section{A Lower Bound on Expansion}
\label{prf}
Let $(X,S,G,\cb)$ be a building-like complex. For a $k$-simplex $\tau=(v_0,\ldots,v_k)$ and $0 \leq i \leq k$, let $\tau_i=(v_0,\ldots,v_{i-1},v_{i+1},\ldots,v_k)$.
The proof of Theorem \ref{gromov} depends on the following homological filling property.
\begin{proposition}
\label{fill}
There exists a family of chains
$$
\ccc=\{c_{s,\tau} \in C_{k+1}(\bst) : -1 \leq k \leq n-1~,~ (s,\tau) \in \cf_k \}
$$
such that
\begin{equation}
\label{cond1}
\partial_{k+1} c_{s,\tau}=\tau+\sum_{i=0}^k c_{s,\tau_i}.
\end{equation}
\end{proposition}
\noindent
{\bf Proof:} We define the $c_{s,\tau}$'s by induction on $k$. First let $k=-1$ and let $*$ be the empty simplex.
For each $s \in S$, choose an arbitrary vertex $v_{s} \in B_{s,*}$
and let $c_{s,*}=v_{s}$. For the induction step, let $0 \leq k \leq n-1$ and suppose that
the $c_{s,\tau}$'s have been defined for all $(s,\tau) \in \cup_{j < k} \cf_j$  and that the family $\{c_{s,\tau}: (s,\tau) \in \cf_j \, , \, -1 \leq j < k\}$ satisfies (\ref{cond1}). Let $(s,\tau) \in \cf_k$. Then
$$z=\tau+ \sum_{i=0}^k c_{s,\tau_i} \in C_k(\bst)+\sum_{i=0}^k C_k(B_{s,\tau_i})\subset C_k(\bst).$$
We claim that $z \in Z_k(\bst)$. Indeed
\begin{equation*}
\begin{split}
\partial_k z &=\partial_k \tau +\sum_{i=0}^k \partial_k c_{s,\tau_i} \\
&=\sum_{i=0}^k \tau_i + \sum_{i=0}^k  (\tau_i+\sum_j c_{s,\tau_{ij}}) \\
&= \sum_{i,j} c_{s,\tau_{ij}} =0.
\end{split}
\end{equation*}
The last equality follows from the fact that each $c_{s,\tau_{ij}}$ appears twice.
As $\thh_k(\bst)=0$, it follows that there exists a $(k+1)$-chain
$c_{s,\tau} \in C_{k+1}(\bst)$  such that $\partial_{k+1} c_{s,\tau} =z$.
It is clear that the family $\{c_{s,\tau}: (s,\tau) \in \cf_j \, , \, -1 \leq j \leq k\}$
satisfies (\ref{cond1}).
{\enp}
\noindent
For $0 \leq k \leq n$ and $s \in S$, define the contraction operator $$\isig :C^{k}(X) \rightarrow
C^{k-1}(X)$$ as follows. For $\alpha \in C^{k}(X)$ and $\tau \in X(k-1)$ let
$$
\isig \alpha (\tau)=\alpha(c_{s,\tau}).
$$
\begin{claim}
\label{homotopy}
For $0 \leq k \leq n-1$ and $\alpha \in C^k(X)$
\begin{equation}
\label{contra}
d_{k-1} \isig \alpha +\isig d_k \alpha =\alpha.
\end{equation}
\end{claim}
\noindent
{\bf Proof:} Let $\tau \in X(k)$. Then
\begin{equation*}
\begin{split}
d_{k-1} \isig \alpha (\tau) +  \isig d_k \alpha(\tau)
&= \isig \alpha (\partial_k \tau) + d_k \alpha (c_{s,\tau}) \\
&= \sum_{i=0}^k \isig \alpha (\tau_i) + \alpha (\partial_{k+1}c_{s,\tau}) \\
&=\sum_{i=0}^k \alpha(c_{s,\tau_i}) + \alpha(\tau+\sum_{i=0}^k c_{s,\tau_i}) \\
&= \sum_{i=0}^k \alpha(c_{s,\tau_i})+\alpha(\tau)+ \sum_{i=0}^k \alpha(c_{s,\tau_i})=\alpha(\tau).
\end{split}
\end{equation*}
{\enp}
\noindent
{\bf Remark:} If $\alpha$ is a $k$-cocycle, then (\ref{contra}) gives a way of representing $\alpha$ as a $k$-cobounday,
i.e. $\alpha=d_{k-1} \isig \alpha$. For a general $\alpha \in C^k(X)$, it provides a another representative of $[\alpha] \in C^k(X)/B^k(X)$.
\ \\ \\
{\bf Proof of Theorem \ref{gromov}.}
Let $0 \leq k \leq n-1$ and $\alpha \in C^k(X)$. Fix $s \in S$ then by Claim \ref{homotopy},
$\isig d_k \alpha= \alpha - d_{k-1} \isig \alpha$. Therefore
\begin{equation}
\label{upper}
\| [\alpha]\| \leq \|\isig d_k \alpha \|.
\end{equation}
For $\eta \in X(k+1)$, let
$$\lambda(\eta)=\frac{1}{|S| \cdot w(\eta)} \sum_{\{ (s,\tau) \in \cf_k \, :\, \eta \in \scst\}} w(\tau)$$
and
$$\tlam(\eta)=\frac{1}{|S| \cdot w(\eta)} \sum_{\{ (s,\tau) \in \cf_k \, :\, \eta \in \bst\}} w(\tau).$$
Let
\begin{equation}
\label{deftheta}
\theta_k=\theta_k(X,\ccc)=\max_{\eta \in X(k+1)} \lambda(\eta).
\end{equation}
\begin{proposition}
\label{lbda}
For $0 \leq k \leq n-1$,
\begin{equation}
\label{lbdae}
h_k(X) \geq \frac{1}{\theta_k}.
\end{equation}
\end{proposition}
\noindent
{\bf Proof:}
Let $\alpha \in C^k(X)$.
Summing (\ref{upper}) over all $s \in S$ we obtain
\begin{equation*}
\begin{split}
|S| \cdot \|[\alpha]\| &\leq \sum_{s \in S} \|\isig d_k \alpha \| \\
&=\sum_{s \in S} \sum \{w(\tau): \tau \in X(k)~,~ \isig d_k \alpha (\tau) \neq 0\} \\
&=\sum_{s \in S} \sum \{ w(\tau): \tau \in X(k)~,~ d_k \alpha (c_{s,\tau}) \neq 0 \}  \\
&\leq \sum_{s \in S} \sum \{w(\tau): \tau \in X(k)~,~ \supp(d_k \alpha) \cap \scst \neq \emptyset \}  \\
&\leq \sum_{\eta \in \supp(d_k \alpha)} \sum_{s \in S}\sum\{w(\tau): \tau \in X(k), \eta \in \scst\} \\
&=\sum_{\eta \in \supp(d_k \alpha)}\sum_{\{ (s,\tau) \in \cf_k \, :\, \eta \in \scst\}} w(\tau) \\
&=\sum_{\eta \in \supp(d_k \alpha)} |S| \cdot w(\eta) \lambda(\eta) \\
&\leq |S|\cdot \theta_k  \sum_{\eta \in \supp(d_k \alpha)} w(\eta) \\
&= |S|\cdot \theta_k \cdot \|d_k \alpha\|.
\end{split}
\end{equation*}
{\enp}
\noindent
To complete the proof of
Theorem \ref{gromov}, it thus suffices to show the following:
\begin{claim}
\label{wtau}
\begin{equation}
\label{bwtau}
\theta_k \leq \binom{n+1}{k+2} a_{k}.
\end{equation}
\end{claim}
\noindent
{\bf Proof:} Fix an $\eta \in X(k+1)$. By the homogeneity condition (C2), $\tlam(\eta)=\tlam(g\eta)$ for all $g \in G$. The transitivity assumption (C1) implies
that $$f_n(X) \leq (G:\stab_G(\eta)) c(\eta)$$
and hence
$$|\stab_G(\eta)| \leq \frac{|G| \cdot c(\eta)}{f_n(X)} =|G|\cdot \binom{n+1}{k+2}w(\eta).$$
Therefore
\begin{equation*}
\begin{split}
|G|\lambda(\eta) \leq |G|\tlam(\eta)&=\sum_{g \in G} \tlam(g\eta) \\
&= \sum_{g \in G}\frac{1}{|S| \cdot w(g\eta)}\sum_{\{(s,\tau) \in \cf_k: g\eta \in \bst\}} w(\tau) \\
&= \frac{1}{|S| \cdot w(\eta)}\sum_{(s,\tau) \in \cf_k}w(\tau) |\{g \in G : g\eta \in \bst\}| \\
&\leq \frac{1}{|S| \cdot w(\eta)}\sum_{(s,\tau) \in \cf_k}w(\tau)\cdot |\stab_G(\eta)| \cdot a_k \\
&=\frac{a_k \cdot |\stab_G(\eta)|}{w(\eta)}\cdot \frac{1}{|S|}\sum_{(s,\tau) \in \cf_k}w(\tau)  \\
&=\frac{a_k \cdot |\stab_G(\eta)|}{w(\eta)} \\
&\leq |G| \cdot \binom{n+1}{k+2} \cdot a_k.
\end{split}
\end{equation*}
{\enp}

\section{Building-Like Complexes}
\label{blike}
In this section we give applications of Theorem \ref{gromov} to two families of building-like complexes.

\subsection{Basis-Transitive Matroidal Complexes}
\label{matroid}
Let $M$ be a matroid on the vertex set $V$ with rank function $\rho$ and let $n=\rho(V)-1$. We identify $M$ with
its $n$-dimensional matroidal complex, namely the simplicial complex on $V$ whose simplices are the
independent sets of the matroid.  Matroidal complexes are characterized by the property that their induced subcomplexes $M[S]$ are pure for every $S \subset V$.
It is well known (see e.g. Theorem 7.8.1 in \cite{B95}) that
$\thh_i(M)=0$ for all $i < \dim M =n$.

A matroid $M$ is {\it basis-transitive} if its automorphism group $Aut(M)$ is transitive on the bases (i.e. maximal faces) of
$M$. One such example is the independence matroid of a vector space. For a classification of basis-transitive matroids see \cite{DL94} and the references therein.
Let $M$ be a basis-transitive matroid of rank $n+1$ and let $G$ be a subgroup of $Aut(M)$ such that $G$ is transitive on the facets. Let $S=M(n)$ be the $G$-set of all $n$-faces of $M$. For $(s,\tau) \in S \times M(k)= \cf_k$ let $\bst=M[s \cup \tau]$.
Then $g\bst=B_{g\sigma,g\tau}$ for all $g \in G$. As $\rho(s \cup \tau)=n+1$, it follows that $\thh_i(\bst)=\thh_i(M[s \cup \tau])=0$ for all $i<n$. Letting
$$\cb=\{\bst: -1 \leq k <n, (s,\tau) \in \cf_k\}$$
it follows that $(M,S,G,\cb)$ is an $n$-dimensional building-like complex.
Now
$$
a_k=a_k(M,S,G,\cb) \leq \max\{f_{k+1}(\bst): (s,\tau) \in \cf_{k}\} \leq \binom{n+k+2}{k+2}.
$$
Writing $\epsilon_1(n,k)=\left(\binom{n+1}{k+2} \binom{n+k+2}{k+2}\right)^{-1}$, Theorem \ref{gromov} implies the following
\begin{corollary}
\label{matexp}
If $M$ is basis-transitive matroid of rank $n+1$ then for all $0 \leq k \leq n-1$,
\begin{equation*}
\label{expmate}
h_k(M) \geq \epsilon_1(n,k).
\end{equation*}
\end{corollary}
\noindent
{\bf Remark}: The bound given in Corollary \ref{matexp} is in general weak and can sometimes be significantly improved for specific classes of basis-transitive matroids
by explicitly constructing a family of chains $\{c_{s,\tau}\}$ satisfying (\ref{cond1}) and then using Proposition \ref{lbda} directly. We illustrate this by the following example.
\ \\ \\
{\bf The Partition Matroid}
\\
Let $V_1,\ldots,V_{n+1}$ be $n+1$ disjoint sets such that $|V_i|= m$ and let $X=X_{n,m}$ be the partition matroid with respect to $V_1,\ldots,V_{n+1}$,
i.e. $\sigma \in X_{n,m}$ iff $|\sigma \cap V_i| \leq 1$ for all $1 \leq i \leq n+1$.
Fix a vector $\uv=(v_1,\ldots,v_{n+1} ) \in \underline{V}=V_1 \times \cdots \times V_{n+1}$.
For an integer $\ell \geq 1$ let $[\ell]=\{1,\ldots,\ell\}$.
Let $-1 \leq k \leq n-1$ and let $\tau=\{u_i:i \in I\} \in X_{n,m}(k)$
where $u_i \in V_i$ and $I \in \binom{[n+1]}{k+1}$.
Define $$j=j(\tau)=\max\{\ell: [\ell] \subset I\}$$
and let
$$\tau'=\{u_i:i \in [j]\} \in X_{n,m}(j-1)~,~\tau''=\{u_i:i \in I-[j]\} \in X_{n,m}(k-j).$$
For $T \subset [j]$, let
$$\sigma_T=\{v_t:t \in T\} \cup\{u_t:t \in [j]-T\}$$
and let $z_{\tau}=\sum_{T \subset [j]}\sigma_T$.
If $v_i \neq u_i$ for all $i \in [j]$, then $z_{\tau}$ is the fundamental cycle of the octahedral  $(j-1)$-sphere
$\{u_1,v_1\}*\cdots* \{u_j,v_j\}$. Otherwise  $z_{\tau}=0$.
Define $\ctil_{\tau} \in C_{k+1}(X_{n,m})$ as the concatination $z_{\tau}v_{j+1}\tau''$.
For $i \in I$, let $\tau_i=\tau-\{u_i\}$.
\begin{claim}
\label{ctau}
$$
\partial_{k+1}\ctil_{\tau}=\tau+\sum_{i \in I} \ctil_{\tau_i}.
$$
\end{claim}
\noindent
{\bf Proof:} Note that for any $i \in [j]$,
$$
\sum_{\{T \subset [j]: \max T=i\}} \sigma_T\tau''=\ctil_{\tau_i}.$$
As $\partial_{j-1}z_{\tau}=0$ it follows that
\begin{equation*}
\begin{split}
\partial_{k+1}\ctil_{\tau}&=\partial_{k+1}(z_{\tau}v_{j+1}\tau'')
=z_{\tau} \partial_{k+1-j}(v_{j+1}\tau'') \\
&=z_{\tau}\tau''+\sum_{i \in I-[j]} z_{\tau}v_{j+1}(\tau''-\{u_i\}) \\
&=z_{\tau}\tau''+\sum_{i \in I-[j]} \ctil_{\tau_i} \\
&=\tau'\tau''+\sum_{\emptyset \neq T \subset [j]} \sigma_T \tau''+\sum_{i \in I-[j]}\ctil_{\tau_i} \\
&=\tau+\sum_{i \in [j]} \left(\sum_{\{T \subset [j]: \max T=i\}} \sigma_T\tau''\right)+\sum_{i \in I-[j]}\ctil_{\tau_i} \\
&= \tau+ \sum_{i \in [j]} \ctil_{\tau_i}+\sum_{i \in I-[j]} \ctil_{\tau_i} \\
&=\tau+\sum_{i \in I} \ctil_{\tau_i}.
\end{split}
\end{equation*}
{\enp}
\noindent
Keeping the notation $j=j(\tau)$, we next note that 
\begin{equation}
\label{nscst}
|\supp(\ctil_{\tau})|=
\left\{
\begin{array}{cl}
2^{j} &  {\rm if~} u_{t} \neq v_t {\rm~for~all~} t \in [j]  \\
0 & {\rm otherwise.}
\end{array}
\right.~~
\end{equation}
Therefore
\begin{equation}
\label{fixs}
\begin{split}
\sum_{\tau \in X_{n,m}(k)} |\supp(\ctil_{\tau})|&=\sum_{j=0}^{k+1} 2^j \cdot \left((m-1)^j\binom{n-j}{k+1-j} m^{k+1-j}\right) \\
&= m^{k+1}\sum_{j=0}^{k+1}\left(\frac{2(m-1)}{m}\right)^j \binom{n-j}{n-k-1}.
\end{split}
\end{equation}
\ \\ \\
Let $S=Aut(X_{n,m})$ be the automorphism group of $X_{n,m}$. For $s \in S$ and $\tau \in X_{n,m}(k)$, let
$c_{s,\tau}=s^{-1} \ctil_{s\tau}$.
Claim \ref{ctau} implies that the family
$$
\ccc=\{c_{s,\tau} \in C_{k+1}(\bst) : -1 \leq k \leq n-1~,~ (s,\tau) \in S \times X_{n,m}(k)\}
$$
satisfies (\ref{cond1}). 
We proceed to compute $\theta_k=\theta_k(X_{n,m},\ccc)$ as defined in (\ref{deftheta}).
First note that $c_{s,s_0\tau}= s_0 c_{ss_0,\tau}$ for all $s_0, s \in S$ and $\tau \in X_{n,m}(k)$.
This, together with the transitivity of $S$ on $X_{n,m}(k)$ and (\ref{fixs}), imply that for all $s \in S$
\begin{equation}
\label{fixss}
\begin{split}
\sum_{\tau \in X_{n,m}(k)} |\supp(c_{s,\tau})| &= \sum_{\tau \in X_{n,m}(k)} |\supp(c_{{\rm identity},\tau})|=\sum_{\tau \in X_{n,m}(k)} |\supp(\ctil_{\tau})| \\
&=  m^{k+1}\sum_{j=0}^{k+1}\left(\frac{2(m-1)}{m}\right)^j \binom{n-j}{n-k-1}.
\end{split}
\end{equation}
The transitivity of $S$ on $X_{n,m}(k+1)$ implies that $\lambda(\eta)$ is independent of $\eta \in X_{n,m}(k+1)$.
As
\begin{equation*}
\label{hmgn}
\begin{split}
\lambda(\eta)&=\frac{1}{|S| \cdot w(\eta)} \sum_{\{ (s,\tau) \in \cf_k \, :\, \eta \in \scst\}} w(\tau) \\
&=\frac{f_{k+1}(X_{n,m})}{|S|f_k(X_{n,m})}|\{(s,\tau) \in \cf_k:\eta \in \supp(c_{s,\tau})\}|,
\end{split}
\end{equation*}
it follows by (\ref{fixs}) that
\begin{equation}
\label{lampart}
\begin{split}
\theta_k &= \frac{1}{f_{k+1}(X_{n,m})} \sum_{\eta \in X_{n,m}(k+1)} \lambda(\eta) \\
&=\frac{1}{|S|f_k(X_{n,m})}\sum_{s \in S} \sum_{\tau \in \cf_k} |\supp(c_{s,\tau})| \\
&=\frac{1}{\binom{n+1}{k+1}}\sum_{j=0}^{k+1} \left(\frac{2(m-1)}{m}\right)^j\binom{n-j}{n-k-1}.
\end{split}
\end{equation}
Proposition \ref{lbda} and (\ref{lampart}) imply the following:
\begin{theorem}
\label{expcolor}
For $0 \leq k \leq n-1$,
\begin{equation}
\label{boundpart}
h_k(X_{n,m}) \geq \frac{\binom{n+1}{k+1}}{\sum_{j=0}^{k+1} (\frac{2(m-1)}{m})^j \binom{n-j}{n-k-1}}.
\end{equation}
\end{theorem}
\noindent
We note some special cases of Theorem \ref{expcolor}.
\\
(i) Let $m=1$. Then $X_{n,1}=\Delta_n$ and
$$h_k(X_{n,1}) \geq \frac{\binom{n+1}{k+1}}{\binom{n}{k+1}}=\frac{n+1}{n-k}$$
thereby recovering the bound (\ref{expsim}).
\\
(ii)
Let $m=2$. Then $X_{n,2}$ is the octahedral $n$-sphere and
$$h_k(X_{n,2}) \geq \frac{\binom{n+1}{k+1}}{\sum_{j=0}^{k+1} \binom{n-j}{n-k-1}}=\frac{\binom{n+1}{k+1}}{\binom{n+1}{n-k}}=1.
$$
This coincides with the result of Proposition 5.5 in \cite{DK12}.
\\
(iii) For general $m$ and $k=n-1$,
$$h_{n-1}(X_{n,m}) \geq \frac{n+1}{\sum_{j=0}^{n} (\frac{2(m-1)}{m})^j}.$$
This is a small improvement over the bound $h_{n-1}(X_{n,m}) \geq \frac{n+1}{2^{n+1}-1}$
given in Proposition 5.7 in \cite{DK12}.
\ \\ \\
We conclude this section with an upper bound on the expansion of $X_{n,m}$.
\begin{claim}
\label{lbpar}
Let $0 \leq k \leq n-1$. If $(k+2)|m$ then
$h_{k}(X_{n,m}) \leq 1$.
\end{claim}
\noindent
{\bf Proof:}
Let $V_i=\bigcup_{j=1}^{k+2} V_{ij}$ where $|V_{ij}|=\frac{m}{k+2}$.
Let $\alpha \in C^{k}(X_{n,m})$ be the indicator function of the following set of $k$-simplices:
$$
\bigcup_{1 \leq i_1< \cdots < i_{k+1} \leq n+1}\bigcup_{\pi \in S_{k+1}} V_{i_1,\pi(1)} \times \cdots \times V_{i_{k+1},\pi(k+1)}.
$$
Then $$|\supp(\alpha)|= \binom{n+1}{k+1} \left(\frac{m}{k+2}\right)^{k+1} (k+1)!.$$ The support of the coboundary of $\alpha$ is
$$B=\supp(d_k \alpha)=\bigcup_{1 \leq i_1< \cdots < i_{k+2} \leq n+1}\bigcup_{\pi \in S_{k+2}} V_{i_1,\pi(1)} \times \cdots \times V_{i_{k+2},\pi(k+2)}$$
and so, $$|B|=\binom{n+1}{k+2} \left(\frac{m}{k+2}\right)^{k+2} (k+2)!.$$
We claim that $\|[\alpha]\|=\|\alpha\|$. Indeed, suppose that $\alpha'=\alpha+d_{k-1}\psi$ where
$\psi \in C^{k-1}(X_{n,m})$. Let
$$C=\{(\sigma,\tau) \in \supp(\alpha') \times B: \sigma \subset \tau\}.$$
As $B=\supp(d_k\alpha)=\supp(d_k\alpha')$, it follows that $|C| \geq |B|$. On the other hand, any $\tau \in X_{n,m}(k)$ is contained in at most $(n-k)\cdot\frac{m}{k+2}$
simplices of $B$. It follows that
$$|B| \leq |C| \leq |\supp(\alpha')| \cdot(n-k)\cdot\frac{m}{k+2}.$$
Therefore
$$|\supp(\alpha')|\geq |B| \frac{k+2}{m(n-k)}=|\supp(\alpha)|.$$
It follows that
\begin{equation*}
\begin{split}
\frac{\|d_k\alpha\|}{\|[\alpha]\|}&=\frac{f_k(X_{n,m})\cdot |\supp(d_k \alpha)|}{f_{k+1}(X_{n,m})\cdot |\supp(\alpha)|} \\
&=\frac{\binom{n+1}{k+1} m^{k+1} \cdot \binom{n+1}{k+2} \left(\frac{m}{k+2}\right)^{k+2} (k+2)!}
{\binom{n+1}{k+2} m^{k+2} \cdot \binom{n+1}{k+1} \left(\frac{m}{k+2}\right)^{k+1} (k+1)!}=1.
\end{split}
\end{equation*}
{\enp}

\subsection{Spherical Buildings}
\label{sbuild}
In this section we use Theorem \ref{bound} to recover Gromov's \cite{G10} uniform lower bound on the expansion of spherical buildings of rank $n+1$. Our notation and terminology follows \cite{B84}.
Let $G=\langle B,N\rangle$ be a finite group with a BN-pair of rank $n+1$ and let $(W,{\mathfrak S})$ be the associated Coxeter system. Here
$W=N/(B \cap N)$ is the Weyl group and ${\mathfrak S}$ is the distinguished set of $n+1$ generators of $W$.  For $J \subset {\mathfrak S}$, let $W_J=\langle J \rangle$
and let $G_J=B W_J B$ be the associated standard parabolic group. For ${\mathfrak s} \in {\mathfrak S}$, let $({\mathfrak s})={\mathfrak S}-\{{\mathfrak s}\}$.
The spherical building $\Delta=\Delta(G;B,N)$ is the $n$-dimensional pure simplicial complex on the vertex set $V=\bigcup_{{\mathfrak s} \in {\mathfrak S}} G/G_{({\mathfrak s})}$ whose
maximal faces, called {\it chambers}, are $C_g=\{gG_{({\mathfrak s})}:{\mathfrak s} \in {\mathfrak S}\}$. Two chambers are {\it adjacent} if their intersection is $(n-1)$-dimensional. For $g \in G$, let $V_g=\{gwG_{({\mathfrak s})}: w \in W, {\mathfrak s} \in {\mathfrak S}\}$. The {\it apartment} $A_g$ is the induced complex
$\Delta[V_g]$. It is a simplicial $n$-sphere whose chambers are $\{C_{gw}\}_{w \in W}$, hence $f_n(A_g)=|W|$.
Any two simplices $\sigma,\tau \in \Delta$ are contained in some apartment $A_g$.
\begin{claim}
\label{intapar}
Let $g_1,\ldots,g_k \in G$ and let $Y=\bigcap_{i=1}^k A_{g_i}$. If $\dim Y=n$ then
$\thh_i(Y)=0$ for all $i \leq n-1$.
\end{claim}
\noindent
{\bf Proof:}
It is convenient to identify  the complex $\Delta$ with its geometric realization. Recall the following:
\begin{itemize}
\item A gallery connecting two simplices $\sigma$ and $\tau$ is a sequence $C_0,C_1,\dots,C_r$ of adjacent chambers so that $\sigma$ is a face of $C_0$ and $\tau$ is a face of $C_r$. The gallery is called minimal if it has minimal length among all possible galleries connecting $\sigma$ and $\tau$.
\item An apartment $A$ which contains two simplices contains also every minimal gallery connecting them.
\item Let $x$ and $y$ be two points in the geometric realization of $\Delta$. Let $A$ be an apartment containing $x$ and $y$.  Consider the sequence of consecutive chambers visited by a minimal geodesic on the sphere $A$ connecting $x$ and $y$. This sequence forms a minimal gallery and hence is contained in any apartment containing both $x$ and $y$.
\end{itemize}
Fix some $x\in Y=\bigcap_{i=1}^k A_{g_i}$. If $Y$ contains a point antipodal to $x$ in some apartment $A_0$ containing $x$ then it follows from the above that $Y$ contains the apartment $A_0$ and hence it follows that $Y=A_0$ which is an $n$-dimensional sphere and the claim holds.
Otherwise it follows that for each $y\in Y$ there is a unique geodesic arc connecting $x$ to $y$ in all the apartments  containing $x$ and $y$ and in particular this geodesic arc is contained in $Y$. This implies that $Y$ is contractible.
{\enp}

\noindent
Let $S=\Delta(n)$ be the set of chambers of $\Delta$.
For $(s,\tau) \in S \times \Delta(k)=\cf_k$, let $\bst=\cap\{A_g:  s, \tau \in A_g\}$.
Letting
$$\cb=\{\bst: -1 \leq k <n, (s,\tau) \in \cf_k\}$$
it follows from Claim \ref{intapar}  that $(\Delta,S,G,\cb)$ satisfies conditions $(C1),(C2)$ and $(C3)$ of Definition \ref{blcom}.
Clearly
$$
a_k=a_k(\Delta,S,G,\cb) \leq f_{k+1}(A_g) \leq \binom{n+1}{k+2} f_n(A_g)= \binom{n+1}{k+2}|W|.
$$
Let $\omega_n$ be the maximal size of a Weyl group of rank $n+1$ and let
$$\epsilon_2(n,k)=\left(\binom{n+1}{k+2}^2 \omega_n\right)^{-1}.$$
Theorem \ref{gromov} then implies
\begin{corollary}
\label{bilexp}
If $G=\langle B,N\rangle$ is a finite group with BN-pair of rank $n+1$, then for all $0 \leq k \leq n-1$
\begin{equation*}
\label{expmate}
h_k(\Delta(G;B,N)) \geq \epsilon_2(n,k).
\end{equation*}
\end{corollary}
\noindent
{\bf Example:} Let $G=GL_{n+2}(\FF_q)=\langle B,N\rangle$ where $B$ is the group of upper diagonal matrices and $N$ is the group of monomial matrices.
The Weyl group of $G$ is the symmetric group $W=\ccs_{n+2}$. The $n$-dimensional spherical building $\Delta=\Delta(G;B,N)$, denoted by
$A_{n+1}(\FF_q)$, is isomorphic to the order complex of all nontrivial linear subspaces of $\FF_q^{n+2}$.
Corollary \ref{bilexp} implies that
for $0 \leq k \leq n-1$,
\begin{equation}
\label{latflat}
h_k\left(A_{n+1}(\FF_q)\right) \geq \left(\binom{n+1}{k+2}^2(n+2)! \right)^{-1}.
\end{equation}
In particular
\begin{equation}
\label{hnmo}
h_{n-1}\left(A_{n+1}(\FF_q)\right) \geq \frac{1}{(n+2)!}.
\end{equation}
\noindent
{\bf Remark:} The uniform lower bound (\ref{latflat}) on the expansion of $A_{n+1}(\FF_q)$ depends on the particular normalization used in the definition of the norm in $C^k(X)$.
Indeed, (\ref{latflat}) fails to hold if the weight of a $k$-simplex is simply taken as $\frac{1}{f_k(X)}$. For example, let $k=0$ and fix an $n$ such that $n+2$ be divisible by $12$.
If $U$ is an $\frac{n+2}{2}$-dimensional subspace of $\FF_q^{n+2}$, then the degree of $U$ in
the underlying graph $G$ of $A_{n+1}(\FF_q)$ is at most
$$2f_0(A_{\frac{n}{2}}(\FF_q)) = q^{\frac{(n+2)^2}{16}(1+o(1))}.$$
This is much smaller than
$$\frac{f_{1}(A_{n+1}(\FF_q))}{f_0(A_{n+1}(\FF_q))}= \frac{q^{\frac{(n+2)^2}{3}(1+o(1))}}{q^{\frac{(n+2)^2}{4}(1+o(1))}}=
q^{\frac{(n+2)^2}{12}(1+o(1))}.$$
It follows that the $1$-dimensional skeleton of $A_{n+1}(\FF_q)$ is not an expander if one uses the normalization giving the same weight to all $i$-simplices.

\section{Applications}
\label{appl}
Lower bounds on coboundary expansion give rise to applications in two directions: topological overlapping and property testing.

\subsection{Topological Overlapping}
\label{s:top}

Let $X$ be a finite $n$-dimensional pure simplicial complex and let $M$ be an $n$-dimensional $\Int_2$-manifold. For a continuous map $f:X \rightarrow M$
and a point $p \in M$, let
$$
\gamma_f(p)=|\{\sigma \in X(n): p \in f(\sigma)\}|.
$$
The following result is due to Gromov \cite{G10}. See also \cite{MW11} for a detailed exposition (including some improved constants) for the case
$M=\Rea^n$ and $X=\Delta_N^{(n)}$.
\begin{theorem}[\cite{G10}]
\label{overlap}
For any $\epsilon>0$ there exists a $\delta=\delta(M,\epsilon)>0$ such that if
$h_k(X) \geq \epsilon$ for all $0 \leq k \leq n-1$, then there exists a point $p \in M$ such that
$\gamma_f(p) \geq \delta f_n(X)$.
\end{theorem}
\noindent
Let $(X,S,G,\cb)$ be an $n$-dimensional building-like complex and let
$$a(X,S,G,\cb)=\max_{0 \leq k \leq n-1} a_k(X,S,G,\cb).$$
The following consequence of Theorem \ref{gromov} was already noted by Gromov (section 2.13 in \cite{G10}) when $X$ is a spherical building
or a partition matroid.
\begin{corollary}
\label{olsb}
For any $0<c$ and an $n$-dimensional $\Int_2$-manifold $M$, there exists a constant $\delta=\delta(c,M)>0$ such that
if $a(X,S,G,\cb)\leq c$, then for any continuous map $f:X \rightarrow M$ there exists a point
$p \in M$ such that $\gamma_f(p) \geq \delta f_n(X)$.
\end{corollary}

\subsection{Property Testing}
\label{s:top}

\begin{definition}
\label{testability}
Let A be a finite set, and let $dist(*,*)$ be a metric on $A^m$. Let $W_m$ a subset of $A^m$
and $P_m$ a subset of $W_m$. Let $\epsilon>0$ and $q \in \Nat$ be fixed.
We say that the membership of $\alpha \in P_m$ (given $\alpha \in  W_m$) is $(q,\epsilon)$-testable,
if there exists a randomized algorithm which queries only q (independent of m) coordinates of $\alpha$ and answers "yes" if $\alpha \in P_m$, while
it answers "no" with probability at least $\epsilon \cdot dist(\alpha,P_m)$.
\end{definition}

In \cite{KL13}, it was observed that coboundary expansion implies that the subspace of coboundaries is testable within the subspace of cochains. The distance function dealt with there was the Hamming distance, but the same applies to the norm used the this paper, provided the algorithm chooses a face with probability equal to its norm. Theorem \ref{gromov}
therefore implies the following.

\begin{corollary}
\label{testc}
For any $0<c$ and $k < n$ there exist an $\epsilon=\epsilon(c,k,n)>0$ such that if an $n$-dimensional building-like complex satisfies
$a_k(X,S,G,\cb) \leq c$, then checking whether a $k$-cochain $\alpha$
is a $k$-coboundary is $(k+2,\epsilon)$-testable.
\end{corollary}

\section{Concluding Remarks}
\label{conc}

We mention some problems related to the results of this paper.
\begin{enumerate}
\item
It would be interesting to improve the bounds given in Theorem \ref{gromov} and its corollaries.
One concrete question is the following. The $1$-dimensional building $A_{2}(\FF_q)$ is the
points vs. lines graph of the Desarguian projective plane of order $q$. It is known that the normalized Cheeger constant of this graph satisfies $h_0(A_{2}(\FF_q))=1-o(1)$ as $q \rightarrow \infty$.
It seems likely that for $n \geq 2$ the bound (\ref{hnmo}) can similarly be improved.
\begin{conjecture}
\label{flagex}
For fixed $n$ and $q \rightarrow \infty$
\begin{equation*}
\label{hnmob}
h_{n-1}\left(A_{n+1}(\FF_q)\right)= 1-o(1).
\end{equation*}
\end{conjecture}

\item
Let $L_n$ be a geometric lattice of rank $n$ with minimal element $\zhat$ and maximal element $\ohat$.
Let $X(L_n)$ be the order complex of $L_n-\{\zhat,\ohat\}$. Then $X(L_n)$ is $(n-2)$-dimensional and
$\thh_k(X(L_n))=0$ for $k \leq n-3$ (see e.g. \cite{B95}).
It would be interesting to find natural families
$\{L_n\}$ for which $h_{n-3}(X(L_n))$ remains uniformly bounded away from zero. For example, is this the case
when $L_n$ is the lattice of partitions of $[n+3]$?

\end{enumerate}
\ \\ \\
{\bf ACKNOWLEDGEMENT}
\\
The authors would like to thank Jake Solomon for helpful discussions.

\end{document}